\documentclass[letterpaper,12pt]{article}
\usepackage{amssymb}
\usepackage{enumitem}
\usepackage[a4paper, left=2.5cm, right=2.5cm, top=2.5cm, bottom=2.5cm, headsep=0.2cm]{geometry} 
\usepackage{epic}
\usepackage{color}

\usepackage{graphicx}
\usepackage{graphics}

\usepackage[latin2]{inputenc}
\newtheorem{df}{Definition}
\newtheorem{lm}{Lemma}
\newtheorem{tw}{Theorem}
\newtheorem{uw}{Remark}
\newtheorem{wn}{Corollary}

\usepackage{amsmath}

\newcommand{\ra}{\rightarrow}

\newcommand{\RA}{\Rightarrow}
\newcommand{\LA}{\Leftrightarrow}

\newcommand{\NN}{\mathbb{N}}

\newcommand{\sq}{\subseteq}
\newcommand{\mm}{\smallsetminus}

\begin{document}

\begin{flushleft}
\textbf{\Large On the distance in some bipartite graphs $L_{k,n}$}\\

$\ $\\
Marcin \L{}azarz\\
Department of Logic and Methodology of Sciences\\
 Wroc\l{}aw University, Poland\\
lazarzmarcin@poczta.onet.pl
\end{flushleft}

$\ $\\
\begin{abstract}
The paper presents some bipartite graph $L_{k,n}$, so called $(k,n)$-level graph, that arise by taking $k$-th and $(n-k)$-th levels of $n$-dimensional Boolean algebra. Two results are establised: (1) precise description of a distance (a shotest path) beteen arbitrary vertices and (2) solution of the problem how many vertices may be reached in $i$ steps starting from some intial point. 
\\
\end{abstract}

$\ $\\


$\ $\\

$\ $\\
\textbf{1. Preliminaries.} For the standard notions such as simple graph, connected graph, path and so see for example  \cite{dist}.
Let us consider a finite and connected simple graph $G=(V,E)$. For arbitrary vertices  $u$ and $v$, let $\left\|uv\right\|$ denote the distance from $u$ to $v$ (the length of a shortest path from $u$ to $v$), and assume that $\left\|uv\right\|=0$ iff $u=v$. 

\begin{df}
\emph{A map $d\colon V\times V\to\NN$ is called a \textbf{metric} in $G$, iff for any vertices $u,v,w$ hold the following conditions:
\begin{eqnarray}
\label{0} d(u,v)=0\ \LA\ u=v,\\
\label{sym} d(u,v)=d(v,u),\\
\label{tro} d(u,v)+d(v,w)\geq d(u,w).
\end{eqnarray}
Then a couple $(G,d)$ is called  a \textbf{metric graph}.} 
\end{df}

\begin{lm}\label{metr}
Let $(G,d)$ be a metric graph, and moreover:
\begin{eqnarray}
\label{1} d(u,v)=1\ \LA\ uv\in E,
\end{eqnarray}
for any $u,v\in V$. Then $d(u,v)\leq\left\|uv\right\|$.
\end{lm}
\textbf{Proof} goes by induction on $n=\left\|uv\right\|$.\\
(1) For $n=0$ or $n=1$ the thesis follows from  (\ref{0}) and (\ref{1}), respectively.\\
(2) For $n>1$ assume induction hypothesis:
$$\forall u,v\in V(\left\|uv\right\|<n\ \RA\ d(u,v)\leq \left\|uv\right\|).$$
(3) Fix $u,v\in V$ such that $\left\|uv\right\|=n$ and take arbitrary path from  $u$ to $v$:
$$u=x_{0}\ra x_{1}\ra\ldots\ra x_{n-1}\ra x_{n}=v.$$
Since $n>1$, there is some intermediate vertex $x_{i}$ (i.e. $0<i<n$), so $\left\|ux_{i}\right\|+\left\|x_{i}v\right\|=\left\|uv\right\|$, because there is no shorter path from $u$ to $v$. Applying the induction hypothesis and  (\ref{tro}) we compute:
$$d(u,v)\leq d(u,x_{i})+d(x_{i},v)\leq\left\|ux_{i}\right\|+\left\|x_{i}v\right\|=\left\|uv\right\|.\ \blacksquare\\$$

\begin{uw}
Condition $(\ref{1})$ is important.
\end{uw}
\textbf{Proof.} Consider a graph $G=(\{1,2,3\},\{12,13,23\})$ and a map $d$ such that: $d(1,2)=d(1,3)=d(2,3)=2.$ It easy to verify that $(G,d)$ is a metric graph. $\blacksquare$\\

$\ $\\
\textbf{2. The distance in a graph $L_{k,n}$.} Let us consider a finite set   $X=\{1,\ldots,n\}$ and fix a natural  $k$ such that  $k<n-k$. Let $[X]^{k}$ stands for the sets of all $k$-element subsets of $X$, and similarly $[X]^{n-k}$ stands for rhe $(n-k)$-elenent subsets of $X$. Let us define a bipartite graph $L_{k,n}=(V,E)$ in the following way:
$$V=[X]^{k}\cup [X]^{n-k},\ \ \ \ E=\{AB:A\in [X]^{k}\wedge B\in [X]^{n-k}\wedge A\subseteq B\}.$$
This graph can be regarded as built-up from two levels of $n$-dimensional Boolean algebra, namely the $k$-th and $(n-k)$-th levels. Hence one can call it a $(k,n)$-\textbf{level graph}.

\begin{tw}
If $A,B\in V$ and $|A\cap B|=i$, then:
\begin{enumerate}
  \item\label{roz} $\left\|AB\right\|\leq 2\left\lceil \frac{k-i}{n-2k}\right\rceil+1$, if $|A|\neq|B|$,
  \item\label{row} $\left\|AB\right\|\leq 2\left\lceil \frac{|A|-i}{n-2k}\right\rceil$, if $|A|=|B|$.
\end{enumerate}
\end{tw}
\textbf{Proof.} For short, put $t=n-2k$. Moreover, $[a,b]$ denotes the set $\{x\in\NN:a\leq x\leq b\}$.

Ad(\ref{roz}). Without loss of the generality assume that $A\in [X]^{k}$, $B\in [X]^{n-k}$, and moreover:
$$A=[1,i]\cup[i+1,k],\ \ \ \ B=[1,i]\cup[k+1,n-i].$$
Define a sequence of sets $C^{0}_{0},C^{0}_{1},C^{1}_{1},\ldots,C^{j}_{j},C^{j}_{j+1},C^{j+1}_{j+1},\ldots,C^{s-1}_{s-1},C^{s-1}_{s},C^{s}_{s}$, in the following way:

$C^{0}_{0}\ \ \ =[1,i]\cup[i+1+0t,k+0t],$

$C^{0}_{1}\ \ \ =[1,i]\cup[i+1+0t,k+1t],$

$\ \ \vdots$

$C^{j}_{j}\ \ \ =[1,i]\cup[i+1+jt,k+jt],$

$C^{j}_{j+1}\,=[1,i]\cup[i+1+jt,k+(j+1)t],$

$\ \ \vdots$

$C^{s}_{s}\ \ \ =[1,i]\cup[i+1+st,k+st],$\\
where $s$ is the smallest natural such that $C^s_s\sq B$, i.e. $s$ satisfies inequalities:
$$k+1\leq i+1+st,\ \ \ \ i+1+(s-1)t<k+1,$$
so $s=\left\lceil\frac{k-i}{t} \right\rceil$. It is easy to observe that $C^{j}_{j}\in [X]^{k}$, $C^{j}_{j+1}\in [X]^{n-k}$ and $C^{j}_{j}\subseteq C^{j}_{j+1}\supseteq C^{j+1}_{j+1}$, and moreover $C^{0}_{0}=A,\ C^{s}_{s}\subseteq B$. Henceforth we just constructed a path in $G$:
$$C^{0}_{0}\ra C^{0}_{1}\ra C^{1}_{1}\ra\ldots\ra C^{s}_{s}\ra B,$$
of the length $2\left\lceil \frac{k-i}{t}\right\rceil+1$, which ends the proof of (\ref{roz}).

Ad(\ref{row}a). Assume that  $A,B\in [X]^{n-k}$ and
$$A=[1,i]\cup[i+1,n-k],\ \ \ \ B=[1,i]\cup[n-k+1,2n-2k-i],$$
and simultaneously  $2n-2k-i\leq n$ i.e. $t\leq i$. Just like before we define sets\\
$C^{0}_{0},C^{1}_{0},C^{1}_{1},\ldots,C^{s-1}_{s-2},C^{s-1}_{s-1},C^{s}_{s-1}$ such that:

$C^{0}_{0}\ \ \ =[1,i]\cup[i+1+0t,n-k+0t],$

$C^{1}_{0}\ \ \ =[1,i]\cup[i+1+1t,n-k+0t],$

$\ \ \vdots$

$C^{j}_{j}\ \ \ =[1,i]\cup[i+1+jt,n-k+jt],$

$C^{j+1}_{j}\,=[1,i]\cup[i+1+(j+1)t,n-k+jt],$

$\ \ \vdots$

$C^{s}_{s-1}\,=[1,i]\cup[i+1+st,n-k+(s-1)t],$\\
where $s$ is the smallest natural such that $C^{s}_{s-1}\sq B$, i.e. $s$ satisfies inequalities:
$$n-k+1\leq i+1+st,\ \ \ \ i+1+(s-1)t<n-k+1,$$
so $s=\left\lceil\frac{n-k-i}{t} \right\rceil$. Then we have: $C^{j}_{j}\in [X]^{n-k}$, $C^{j+1}_{j}\in [X]^{k}$, $C^{j}_{j}\supseteq C^{j+1}_{j}\subseteq C^{j+1}_{j+1}$ and $C^{0}_{0}=A,\ C^{s}_{s-1}\subseteq B$. Finally, exists a path in $G$ from $A$ to $B$ of the length  $2\left\lceil \frac{n-k-i}{t}\right\rceil=2\left\lceil \frac{|A|-i}{t}\right\rceil$:
$$C^{0}_{0}\ra C^{1}_{0}\ra C^{1}_{1}\ra\ldots\ra C^{s}_{s-1}\ra B.$$

Ad(\ref{row}b). Assume that $A,B\in [X]^{k}$ and
$$A=[1,i]\cup[i+1,k],\ \ \ \ B=[1,i]\cup[k+1,2k-i].$$
We will reduce this case to (\ref{roz}). Consider two posibilities: first, if  $k-i\geq t$ then putting $C=[1,i+t]\cup[k+1,2k-i]$, we get $C\in [X]^{n-k}$ and $B\subseteq C$, and moreover $|A\cap C|=i+t$. By part  (\ref{roz}) we have a path from  $A$ to $C$ of the length  $2\left\lceil \frac{k-(i+t)}{t}\right\rceil+1$. Obverve also that:
$$2\left\lceil \frac{k-(i+t)}{t}\right\rceil+1=2\left\lceil \frac{k-i}{t}-1\right\rceil+1=2\left\lceil \frac{k-i}{t}\right\rceil-1;$$
hence we achieve a path from $A$ to $B$ of the length $2\left\lceil \frac{k-i}{t}\right\rceil$ which ends the proof. The second possibility $k-i<t$ is trivial, because then we obtain $|A\cup B|<n-k$ so the path from  $A$ to $B$ is of the length $2=2\cdot1=2\left\lceil \frac{k-i}{t}\right\rceil$.  $\blacksquare$\\

The above proof give an algorithm of constructing a path from arbitrary vertex $A$ to arbitrary $B$. Let us define the function $d$ for $A,B\in V$ (just like before $|A\cap B|=i$):

\begin{displaymath}
d(A,B)=\left\{ \begin{array}{cc}
2\left\lceil \frac{k-i}{n-2k}\right\rceil+1, & \textrm{if $\ |A|\neq|B|$}\\
2\left\lceil \frac{|A|-i}{n-2k}\right\rceil,\ \ \ \ & \textrm{if $\ |A|=|B|$}
\end{array} \right. .
\end{displaymath}

\begin{lm}\label{troo}
The map $d$ is a metric and satisfies $(\ref{1})$.
\end{lm}
\textbf{Proof.} It is easy to see that  $d$ fulfils conditions (\ref{0}), (\ref{sym}) and (\ref{1}). The proof of (\ref{tro}) is also simple but quite arduous (eight cases). Assume abreviations $|(A\cap B)\mm C|=p$, $|(A\cap C)\mm B|=q$, $|(B\cap C)\mm A|=r$, $|A\cap B\cap C|=x$, and moreover $|A\cap B|=i$, $|B\cap C|=j$, $|A\cap C|=l$.\\
(1) Consider the the case when $A,B,C\in [X]^{k}$. To prove $d(A,B)+d(B,C)\geq d(A,C)$ it is sufficient to show that:
$$\frac{k-i}{n-2k}+\frac{k-j}{n-2k}\geq \frac{k-l}{n-2k}\ .$$
However it is clear that the above inequality is equal to:
\begin{equation}
\label{yyy} k+q\geq p+x+r,
\end{equation}
which is obviously true.\\
(2) If $A,C\in [X]^{n-k}$, $B\in [X]^{k}$, it is sufficient to show that:
$$\frac{k-i}{n-2k}+\frac{k-j}{n-2k}+1\geq \frac{n-k-l}{n-2k},$$
which is also equivalent to (\ref{yyy}).\\
(3) Next six cases we easy check in similar way. $\blacksquare$\\

The main result of this section is 

\begin{wn}\label{met}
$d(A,B)=\left\|AB\right\|$, for any $A,B\in V$.
\end{wn}
\textbf{Proof.} Inequality $\leq$ follows forom lemmas \ref{metr} and \ref{troo}; inequality $\geq$ is obvious, since number  $d(A,B)$ is length of concrete path from $A$ to $B$. $\blacksquare$\\

$\ $\\
\textbf{3. The number of vertices reachable in $i$ steps.} The set  $P=\{1,\ldots,k\}$ is called an \textbf{initial vertex} of a graph $L_{k,n}$. Our aim is to find a pattern of the function $f$ that describe a cardinality of the set of vertices, that we reach in consecutive steps, starting from the initial vertex.

To ilustrate the problem let us consider the graph $L_{2,5}$. For simplify notation, vertex $\{a,b\}$ will be denoted $ab$ and similarly $abc$ stands for the vertex $\{a,b,c\}$. 

The initial vertex is the only one that is reached in $0$ steps, so $\Gamma(0)=\{12\}$. In one step we reach three vertices: $123$, $124$, $125$, so we write $\Delta(0)=\{123,124,125\}$. In two steps we reach six vertices: $\Gamma(1)=\{13,14,15,23,24,25\}$ and so on. So the function that we are looking for, in the case of graph $L_{2,5}$ is: $f(0)=1$, $f(1)=3$, $f(2)=6$, $f(3)=6$, $f(4)=3$, $f(5)=1$ (see figure below).

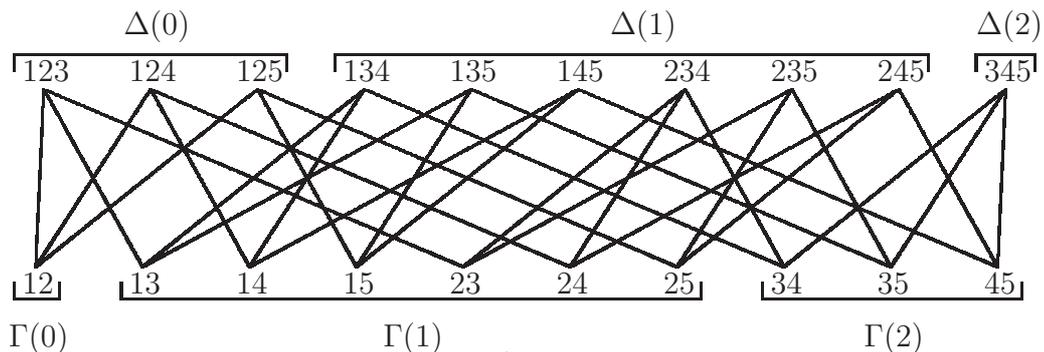
\begin{center}
\begin{picture}(360,110)
\thicklines

\put(0,80){123}   \put(40,80){124}   \put(80,80){125}   \put(120,80){134} \put(160,80){135} \put(200,80){145}
\put(240,80){234} \put(280,80){235}  \put(320,80){245}  \put(360,80){345}

\put(0,0){12} \put(40,0){13} \put(80,0){14} \put(120,0){15} \put(160,0){23} \put(200,0){24} \put(240,0){25} \put(280,0){34} \put(320,0){35} \put(360,0){45}

\qbezier(5,10)(5,10)(8,77)     \qbezier(5,10)(5,10)(48,77)     \qbezier(5,10)(5,10)(88,77)
\qbezier(45,10)(45,10)(8,77)   \qbezier(45,10)(45,10)(128,77)  \qbezier(45,10)(45,10)(168,77)
\qbezier(85,10)(85,10)(48,77)  \qbezier(85,10)(85,10)(128,77)  \qbezier(85,10)(85,10)(208,77)
\qbezier(125,10)(125,10)(88,77)  \qbezier(125,10)(125,10)(168,77)  \qbezier(125,10)(125,10)(208,77)
\qbezier(165,10)(165,10)(8,77)  \qbezier(165,10)(165,10)(248,77)  \qbezier(165,10)(165,10)(288,77)
\qbezier(205,10)(205,10)(48,77)  \qbezier(205,10)(205,10)(248,77)  \qbezier(205,10)(205,10)(328,77)
\qbezier(245,10)(245,10)(88,77)  \qbezier(245,10)(245,10)(288,77)  \qbezier(245,10)(245,10)(328,77)
\qbezier(285,10)(285,10)(128,77)  \qbezier(285,10)(285,10)(248,77)  \qbezier(285,10)(285,10)(368,77)
\qbezier(325,10)(325,10)(168,77)  \qbezier(325,10)(325,10)(288,77)  \qbezier(325,10)(325,10)(368,77)
\qbezier(365,10)(365,10)(208,77)  \qbezier(365,10)(365,10)(328,77)  \qbezier(365,10)(365,10)(368,77)

\qbezier(-3,-2)(-3,-2)(14,-2)     \qbezier(-3,-2)(-3,-2)(-3,4)     \qbezier(14,-2)(14,-2)(14,4)
\qbezier(37,-2)(37,-2)(254,-2)    \qbezier(37,-2)(37,-2)(37,4)     \qbezier(254,-2)(254,-2)(254,4)
\qbezier(277,-2)(277,-2)(374,-2)  \qbezier(277,-2)(277,-2)(277,4)  \qbezier(374,-2)(374,-2)(374,4)
\put(-5,-20){$\Gamma(0)$}          \put(135,-20){$\Gamma(1)$}         \put(315,-20){$\Gamma(2)$}

\qbezier(-3,90)(-3,90)(99,90)     \qbezier(-3,90)(-3,90)(-3,84)     \qbezier(99,90)(99,90)(99,84)
\qbezier(117,90)(117,90)(339,90)     \qbezier(117,90)(117,90)(117,84)     \qbezier(339,90)(339,90)(339,84)
\qbezier(357,90)(357,90)(379,90)     \qbezier(357,90)(357,90)(357,84)     \qbezier(379,90)(379,90)(379,84)
\put(38,98){$\Delta(0)$}          \put(220,98){$\Delta(1)$}         \put(357,98){$\Delta(2)$}
\end{picture}
\end{center}

\vspace{0.9cm}

We assume that the Newton's symbol ${n\choose k}$ have a sense for any integers $n$ and $k$:

\begin{displaymath}
{n \choose k}=\left\{ \begin{array}{cc}
\frac{n!}{k!(n-k)!}, & \textrm{if $\ k\geq0\wedge k\leq n$}\\
0, & \textrm{if $\ k<0\vee k> n$}
\end{array} \right..
\end{displaymath}
Fix  $n$ and $k$ such that $k<n-k$ and assume $t=n-2k$, $s=\left\lceil \frac{k}{t}\right\rceil$. For $i=0,\ldots,s$ put:
$$\Gamma(i)=\{A\in [X]^{k}:d(P,A)=2i\},$$
$$\Delta(i)=\{B\in [X]^{n-k}:d(P,B)=2i+1\}.$$
The set $\Gamma(i)$ is just a set of vertices that may be reached from  $P$, in precisely  $2i$ steps. Similarly, $\Delta(i)$ is a set of vertices that may be reached from  $P$, in precisely  $2i+1$ steps. By Corollary \ref{met} it easy follows that:
$$\Gamma(i)\cap\Gamma(j)=\emptyset,\ \ \ \Delta(i)\cap\Delta(j)=\emptyset, \ \ \emph{\emph{ for }}\ i\neq j,$$
and
$$[X]^{k}=\bigcup^{s}_{i=0}\Gamma(i),\ \ \ [X]^{n-k}=\bigcup^{s}_{i=0}\Delta(i).$$
Let us compute the cardinality of $\Gamma(i)$; first observe that:
$$d(P,A)=2i\ \LA\ 2\left\lceil\frac{k-|P\cap A|}{t} \right\rceil=2i\ \LA\ \left\lceil\frac{|P\smallsetminus A|}{t} \right\rceil=i\ \LA$$
$$i-1<\frac{|P\smallsetminus A|}{t}\leq i\ \LA\ (i-1)t<|P\smallsetminus A|\leq it\ \LA\ \bigvee ^{t}_{j=1}|P\smallsetminus A|=(i-1)t+j,$$
so:
$$\Gamma(i)=\bigcup^{t}_{j=1}\{A\in [X]^{k}: |P\smallsetminus A|=(i-1)t+j\}.$$
The set $\Gamma(i)$ has been presented as a union of  disjoint sets. It is also clear that:
$$|\{A\in [X]^{k}: |P\smallsetminus A|=l\}|={k\choose k-l}{n-k\choose l},$$
so we achieve:
$$\gamma(i)=|\Gamma(i)|=\sum^{t}_{j=1}{k\choose k-((i-1)t+j)}{n-k\choose (i-1)t+j}.$$

Now let us compute the cardinality of  $\Delta(i)$. Similarly like previously we show that:
$$\Delta(i)=\bigcup^{t}_{j=1}\{B\in [X]^{n-k}: |P\smallsetminus B|=(i-1)t+j\},$$
and since
$$|\{B\in [X]^{n-k}: |P\smallsetminus B|=l\}|={k\choose k-l}{n-k\choose t+l},$$
we obtain:
$$\delta(i)=|\Delta(i)|=\sum^{t}_{j=1}{k\choose k-((i-1)t+j)}{n-k\choose it+j}.$$

The main result of this section is

\begin{wn} The function $f$ that gives the cardinality of the set of all vertices that is reached in precisely $x\in\{0,1,\ldots, 2\left\lceil \frac{k}{t}\right\rceil+1\}$ steps is
\begin{displaymath}
f(x)=\left\{ \begin{array}{cc}
\gamma(i), & \textrm{if $\ x=2i\ \ \ \ \ $}\\
\delta(i), & \textrm{if $\ x=2i+1$}
\end{array} \right. .
\end{displaymath}
\end{wn}

From the above calculation we obtain also a pure combinatorial corollary:

\begin{wn}\label{komb} For   $n$ and $k$ such that $2k<n$ hold:
\begin{enumerate}
	\item $\displaystyle{n\choose k}=\sum^{\left\lceil \frac{k}{n-2k}\right\rceil}_{i=0}\sum^{n-2k}_{j=1}{k\choose k-((i-1)(n-2k)+j)}{n-k\choose (i-1)(n-2k)+j},$
	\item $\displaystyle{n\choose k}=\sum^{\left\lceil \frac{k}{n-2k}\right\rceil}_{i=0}\sum^{n-2k}_{j=1}{k\choose k-((i-1)(n-2k)+j)}{n-k\choose i(n-2k)+j}.$\\
\end{enumerate}
\end{wn}

\end{document}